\documentclass{article}
\title{Analysis and Extension of Omega-Rule}
\author{Ryota Akiyoshi, Grigori Mints}
\date{March 2, 2011}
\usepackage{proof}
\usepackage{amscd}
\usepackage{amsmath,amssymb}
\usepackage{mathptmx}
\newtheorem{definition}{Definition}
\newtheorem{theorem}{Theorem}
\newtheorem{lemma}{Lemma}

\newtheorem{proposition}{Proposition}

\newtheorem{remark}{Remark}

\newtheorem{corollary}{Corollary}

\newcommand{\R}{\mathrm{R}}
\newcommand{\E}{\mathrm{E}}
\newcommand{\D}{\mathrm{D}}
\newcommand{\Srm}{\mathrm{S}}
\newcommand{\Rep}{\mathrm{Rep}}

\begin{document}
\maketitle

\begin{abstract} $\Omega$-rule was introduced by W. Buchholz to give an ordinal-free cut-elimination proof for a subsystem of analysis with $\Pi^{1}_{1}$-comprehension. His proof provides cut-free derivations by familiar rules only for arithmetical sequents. When second-order quantifiers are present, they are introduced by $\Omega$-rule and some residual cuts are not eliminated. Using an extension of $\Omega$-rule we obtain (by the same method as W. Buchholz) complete cut-elimination: any derivation of arbitrary sequent is transformed into its cut-free derivation by the standard rules (with induction replaced by $\omega$-rule).

W. Buchholz used $\Omega$-rule to explain how reductions of finite derivations (used by G. Takeuti for subsystems of analysis) are generated by cut-elimination steps applied to derivations with $\Omega$-rule. We show that the same steps generate standard cut-reduction steps for infinitary derivations with familiar standard rules for second-order quantifiers. This provides an analysis of $\Omega$-rule in terms of standard rules and ordinal-free cut-elimination proof for the system with the standard rules for second-order quantifiers. In fact we treat the subsystem of $\Pi^{1}_{1}$-CA (of the same strength as $ID_{1}$) that W. Buchholz used for his explanation of finite reductions. Extension to full $\Pi^{1}_{1}$-CA is forthcoming in another paper.
\end{abstract}

\section{Introduction}

Our main goal is to extend $\Omega$-rule from \cite{Buchholz01} to get cut-elimination to arbitrary (not only arithmetical) end-formulas. As a warm-up we adapt in the Section 3 the proof from \cite{Buchholz01} to obtain cut-elimination by (what we call) standard reductions in an infinitary system with $\omega$-rule.
As a first step we treat in this paper a system
$\textrm{BI}^{-}$ of the strength of $ID_1$ with $\omega$-rule. An extension to
$\Pi^1_1$-analysis with Bar Induction and $\omega$-rule is planned for another paper. 
We use terminology from \cite{Buchholz01} and rely on the results of this paper. 

The first constructive proof of the cut-elimination theorem for $\Pi^{1}_{1}$-analysis has
been given by G. Takeuti \cite{Takeuti67, Takeuti87} who used a new kind of ordinal notations
(ordinal diagrams) to prove the termination of his cut-elimination
steps. Since then several other proofs appeared
\cite{Arai85, Buchholz81a, Buchholz_Schutte88, Howard72,Pohlers81}, but these
proofs also use some kind of complicated ordinal notations. Much more general
proof using computability predicate introduced by J-Y. Girard \cite{Girard87} 
employs much stronger tools. Moreover, the proof by
G. Takeuti generalizes Gentzen's second consistency proof \cite{Gentzen38}
which 
sacrifices transparency to retain finitistic framework. 

Cut-elimination for $\Pi^{1}_{1}$--analysis with $\omega$-rule has been proved
by M. Yasugi \cite{Yasugi70}. She used Takeuti reductions and assignment of
ordinal diagrams for derivations for proving existence of a cut-free
normal form.

The set of reductions  we employ in normalization
procedure for $\textrm{BI}^{-}$ consists of familiar cut-elimination steps
for the 
second-order arithmetic with $\omega$-rule.

In particular the following the derivation $d$:
$$
\infer[Cut_{C}]{\Gamma}{\infer[\bigwedge_{\forall X C_{0}(X)}]{\Gamma,
\forall X C_{0}(X)}{\deduce{\Gamma, C_{0}(X),\forall X
C_{0}(X)}{\vdots}} & \infer[\bigvee^{T}_{\exists X \neg
C_{0}(X)}]{\Gamma, \exists X \neg C_{0}(X)}{\deduce{\Gamma, \neg
C_{0}(T), \exists X \neg C_{0}(X)}{\vdots}}} 
$$
where $C=\forall X C_{0}(X)$
reduces to the following derivation: 






\footnotesize
$$
\infer[Cut_{C}]{\Gamma}{\infer[\bigwedge_{\forall X C_{0}(X)}]{\Gamma,
\forall X C_{0}(X)}{\deduce{\Gamma, C_{0}(X),\forall X
C_{0}(X)}{\vdots}} & \infer[Cut_{C_{0}(T)}]{\Gamma, \exists X \neg
C_{0}(X)}{\infer[\mathcal{S}^{X}_{T}]{\Gamma, C_{0}(T)}{\infer[Cut_{C}]{\Gamma,
C_{0}(X)}{\deduce{\Gamma, C_{0}(X), \forall X C_{0}(X)}{\vdots} &
\deduce{\Gamma, \exists X \neg C_{0}(X)}{\vdots}}} &
\deduce{\Gamma,  \neg C_{0}(T), \exists X \neg C_{0}(X)}{\vdots}} } 
$$

\normalsize
where $\mathcal{S}^X_T$ is the result of substituting $X$ by $T$. 
Here the second-order cut on the formula $\exists X C_0(X)$ derived from
$C_0(T)$ is replaced by a cut on $C_0(X)$ followed by the substitution
rule
with ``parametric'' occurrences of $C,\neg C$ cut out in a familiar way. 

The {\em standard reductions} are {\em permutation} of a cut with an adjacent logical rule, and 
{\em essential reduction} of cut when both cut formulas are introduced by logical rules or when one of the premises is an axiom. 
We reduce 
cut-elimination 
to a similar question for another normalization procedure
devised by W. Buchholz \cite{Buchholz01} for a different axiomatization of $\textrm{BI}^{-}$
using $\Omega$-rule 
for which he gave an ordinal-free proof of the cut-elimination.
The principal tool is the use of his rules
$\E,\D,\Srm^{X}_{T}$  (going back to the second author's \cite{Mints75a}) 
that provide a fine structure and 
allow to give  detailed analysis of cut-elimination process by standard reductions. 
W. Buchholz has shown how normalization of $\Omega$-derivations induces
normalization steps used in Gentzen's second consistency proof for PA
\cite{Gentzen38} and in Takeuti's proof for $\Pi^1_1$-analysis. We establish
similar result for standard reductions closely following Buchholz's schema. 

Let's recall the definition of 
reduction for infinitary derivations \cite{Tait65, Schwichtenberg77}.
\[ d \mapsto d' \ \textrm{if $d'$ is obtained from $d$ by a standard reduction. }\]
$$
\infer[\textrm{for every inference rule} \ \mathcal{I}]{\mathcal{I} \{  d_{i} \} \mapsto  \mathcal{I} \{  d_{i}'\} }{ d_{i} \mapsto d'_{i} \ \textrm{for all relevant} \ i}
$$
We restate the reduction relation mentioned in \cite{Buchholz01} (the end of Section 4) in the form of a normalization tree $T_{d}$ defined for each derivation $d$ in $\textrm{BI}^{-}$. Nodes of this tree are labelled by derivations in $\textrm{BI}^{-}$ in (approximately) the following way. 
\begin{enumerate}

\item $d$ is placed to the root.

\item If $d$ ends in a logical inference (including $\omega$-rule), then the
predecessors of the root are derivations of the premises.

\item If the end-piece of $d$ contains a cut, then the predecessor is the result  of a standard reduction of this cut. 
\end{enumerate}
We call a derivation $d$ in $\textrm{BI}^{-}$ {\em normalizable} iff $T_{d}$ is 
well-founded. This is an obvious analog of the normalization
in a finite number of steps, cf. the definition of {\em Reduzierforschrift} in
\cite{Gentzen36}. The main result of the first part of this paper is 
\begin{description}
\item[Theorem.]  Every derivation in $\textrm{BI}^{-}$ of an arithmetical sequent is normalizable by standard reductions.
\end{description}
(cf. Theorem 2 in Section 3.2).

It is proved by
establishing the well-foundedness of another normalization tree $T_{d}^{+}$ for
derivations $d$ in a wider system $\textrm{BI}$ using ``look-ahead" functions $tp(d)$ and $d[i]$ from \cite{Buchholz01},
and 
detailed comparison of $T^{+}_d$ and $T_{d}$ using finite structure present in
$T^{+}_{d}$. 

In the second part of the paper, we extend the $\Omega$-rule to achieve normalization for derivations of arbitrary sequents in our language. The $\Omega$-rule can be stated as follows:
\begin{equation}
\label{Omega}
\infer{\Gamma, \exists X A(X)}{\infer={\dots d_{q}: \Delta, \Gamma, \exists X A(X) \dots }{q: \Delta, \neg A(X) }} 
\end{equation}
with a separate premise $d_{q}$ for every cut-free derivation of an arithmetical sequent. Now we change it to have a premise $d_{q}$ for every cut-free derivation $q$ of arbitrary sequent. As W. Buchholz pointed out in private correspondence, such change taken literally would interfere with the translation of ordinary second-order existential rule $A(T) / \exists X A(X)$ into $\Omega$-rule creating an infinite loop. To resolve this we use distinction between explicit and implicit formulas introduced by G. Takeuti \cite{Takeuti87}. Implicit formulas in a derivation are those treceable to cut-formulas. Explicit formulas are those traceable to the end-sequent. For example, a derivation of empty sequent cannot contain explicit formulas. 

 We modify the translation $()^{\infty}$ used in \cite{Buchholz01} in the following way. Only implicit $\exists X$-rules are translated into $\Omega$-rules while explicit rules $A(T)/\exists X A(X)$ are left intact. Derivations $q:\Delta, \neg A(X)$ in (\ref{Omega}) are now ``explicit" derivations: all formulas in $\Delta$ should be explicit, and for arithmetical formula $\neg A(X)$ the explicit/implicit distinction is unimportant.

After this change the cut-elimination proof from \cite{Buchholz01} goes through (with suitable additions) for arbitrary (not only arithmetical) end-sequents. To make induction on derivations possible, we have to consider subderivations where some formulas in the end-sequent are to be treated as implicit since 
they are traceable to a cut below the end-sequent. Such a situation is accounted by introduction of marks. Every formula in every sequent is marker by $e$ (explicit) or $i$ (implicit). Inference rules are stated respecting these marks, so that almost every rule has two version. For example,
$$
\infer{\Gamma, (A \vee B)^{e}}{\Gamma, A^{e}, (A \vee B)^{e}}
\qquad
\infer{\Gamma, (A \vee B)^{i}}{\Gamma, A^{i}, (A \vee B)^{i}}
$$

The present paper consists of 4 sections.  In Section 2, we recall basic definitions and Buchholz's  infinitary systems $\textrm{BI}^{\Omega}_{0}$ and $\textrm{BI}^{\Omega}$, and introduce our target systems $\textrm{BI}^{- }$, (which is a parameter-free subsystem of $\Pi^{1}_{1}$-analysis with $\omega$-rule) and $\textrm{BI}$, which is obtained by adding the rules $\E,\D, \Srm^{X}_{T}$. Moreover, we recall some tools like $tp(d), d[i]$ due to Buchholz \cite{Buchholz01}. 

In Section 3 we give an ordinal-free proof of the cut-elimination theorem for $\textrm{BI}^{-}$ and analyse our reductions. We define reduction relation $red$, normalization tree
$T^{+}_d$ for $\textrm{BI}$-derivations $d$, then  normalization tree 
$T_{d}$ for $\textrm{BI}^{-}$-derivations $d$  as a result of deleting (almost all)
fine structure operations. The well-foundedness of $T_{d}$ is proved for derivations of arithmetical sequents in Section 3.1, then the structure of $T_{d}$ is analysed (Section 3.2).

We extend Buchholz's $\Omega$-rule in Section 4 so that the cut-elimination theorem can be proved for derivations of arbitrary sequent. The system $\textrm{BI}^{\Omega^{+}}$ with the extended $\Omega$-rule is introduced based on language with marks $e,i$, and the cut-elimination theorem for $\textrm{BI}^{\Omega^{+}}$ is proved in Section 4.1. After defining the embedding function from $\textrm{BI}^{-}$ to $\textrm{BI}^{\Omega^{+}}$, we prove that any $d \in  \textrm{BI}^{-} $ is translated into a cut-free derivation $d' \in \textrm{BI}^{\Omega^{+}}_{0}$ (Section 4.2).

We acknowledge the help of three anonymous referees who criticized the first version of  this paper.

\section{Preliminaries}

We adopt Buchholz's formulation of a formal language $L$ although free number variables are redundant because all formal systems in this paper contain $\omega$-rule.
The function symbols are $0$ and $S$. We assume the $n$-ary predicate symbol $R$ for an $n$-ary primitive recursive relation. 
\textit{Atomic formula} is of the form $R(t_{1},\dots,t_{n})$ or $X(t)$ where $X$ is a free predicate variable. $A$ and $\neg A$ where $A$ is atomic are called \textit{literals}.
Formulas are obtained from literals by $\wedge, \vee, \exists x, \forall x, \forall x, \exists X, \forall X$ with the restriction 
that $\forall X A$ and $\exists X A$ are formulas {\em only if  $A$ contains no second-order quantifier and no free predicate variable other than $X$.}

If $A$ is a formula which is not atomic, then its \textit{negation} $\neg A$ is defined using De Morgan's laws. The set of literals without free variables true in the standard model is denoted by TRUE.
Formulas which do not contain second-order quantifiers are called \textit{arithmetical}. 
$FV(A)$ denotes the set of free variables of a formula $A$.

$rk(A)$ (the logical complexity of $A$) is defined as follows.

$rk(A):= 0$ if $A$ is a literal , $\forall X B(X)$ or $\exists X B(X)$. 

$rk(A \wedge B) := rk(A \vee B) =  max(rk(A),  rk(B)) + 1$.

$rk(\forall x A(x)) := rk(\exists x A(x)) = rk(A(0)) +1 $.

In inference rules only the \textit{minor formulas} (which occur in the premises of the rule) and the \textit{principal formulas} (which occur in the conclusion of the rule) are shown explicitly. 
Let $I$ be an inference symbol of a system. Then we write $\Delta(I)$ and $|I|$ in order to indicate the set of principal formulas of $I$ and the index set of $I$ respectively. $\bigcup_{i \in |I|}(\Delta_{i}(I))$ is the set of the minor formulas of $I$.
If $d = I(d_{i})_{i \in|I|}$, then $d_{i}$ denotes the subderivation of $d$ indexed by $i$. If $d$ is a derivation, $\Gamma(d)$ denotes its last sequent.  

We use the systems $\textrm{BI}^{\Omega}_{0}$, $\textrm{BI}^{\Omega}$ introduced in \cite{Buchholz01}.  $\textrm{BI}^{\Omega}_{0}$ is just arithmetic with $\omega$-rule and Repetition rule (due to the second author).
$\textrm{BI}^{\Omega}$ is obtained by adding the rules $\Omega_{\neg \forall X A}$ and $\widetilde{\Omega}_{\neg \forall X A}$ to $\textrm{BI}^{\Omega}_{0}$.

The cut-degree $dg(d)$  is defined just as in \cite{Buchholz01}.
Let $d$ be a derivation in $\textrm{BI}^{\Omega}$. As usual, we write 
\[d \vdash_{m} \Gamma
\] 
if the end-sequent $\Gamma(d) \subseteq \Gamma$ and $dg(d) \leq m$. In what follows, we assume that $\Gamma(d) = \Gamma$ unless otherwise noted. 

The operators $\mathcal{R}_{C}, \mathcal{E}, \mathcal{D}, \mathcal{S}^{X}_{T}$ are defined as in \cite{Buchholz01}. We use their properties stated in Theorems 1, 2, 3 of \cite{Buchholz01}.

\subsection{The Systems $\textrm{BI}^{-}$ and  $\textrm{BI}$}

The system $\textrm{BI}^{-}$ is a parameter-free subsystem of $\Pi^1_1$-analysis with $\omega$-rule. Replacing $\omega$-rule with induction axiom would give a system of strength of $ID_{1}$(, which is $\textrm{BI}_{1}^{-}$ in \cite{Buchholz01}). The system $\textrm{BI}$ is obtained by adding the rules $\E, \D, \Srm^{X}_{T}$ to $\textrm{BI}^{-}$. Thus the system obtained from our $\textrm{BI}$ by replacing $\omega$-rule with induction is $\textrm{BI}_{1}^{*}$ in \cite{Buchholz01}.

\begin{definition}
The systems $\textrm{BI}^{-}$ and $\textrm{BI}$
\end{definition}
\begin{enumerate}

\item $\textrm{BI}^{-}$ consists of the following rules:

$$
\deduce[(\textrm{Ax}_{\Delta})]{}{}\ \infer{\Delta}{} \ \textrm{where} \ \Delta = \{ A \} \subseteq \ \textrm{TRUE or } \ \Delta = \{ C, \neg C \}
$$

$$ 
\deduce[\Huge{(\bigwedge_{A_{0} \wedge A_{1}})}]{}{}\ \infer{A_{0} \wedge A_{1}}{A_{0} & A_{1}}
\qquad
\deduce[\Huge{(\bigvee^{k}_{A_{0} \vee A_{1}})}]{}{}\ \infer[\textrm{where} \ k  \in \{ 0, 1 \}]{A_{0} \vee A_{1}}{A_{k}} 
$$

$$
\deduce[\Huge{(\bigwedge_{\forall x A}})]{}{}\ \infer{\forall x A}{\dots A(x/n) \dots \ \textrm{for all} \ n \in \omega}  
\qquad
\deduce[\Huge{(\bigwedge^{k}_{\exists x A})}]{}{}\ \infer[\textrm{where} \ k  \in \omega]{\exists x A}{A(x/k)}
$$

$$
\deduce[\Huge{(\bigwedge_{\forall X A}})]{}{}\ \infer[\textrm{where} \ Y \ \textrm{is an eigenvariable}]{\forall X A}{A(X/Y) } 
\qquad
\deduce[\Huge{(\bigvee^{T}_{\neg \forall X A})}]{}{}\ \infer{\neg \forall X A}{\neg A(X/T)} 
$$

$$
\deduce[(\R_{C})]{}{}\ \infer{\emptyset}{C & \neg C}
$$

\item $\textrm{BI}$ is obtained by adding the following rules to $\textrm{BI}^{-}$.
 
$$
\deduce[(\E)]{}{}\infer{\emptyset}{\emptyset}
\qquad
\deduce[(\D)]{}{}\infer{\emptyset}{\emptyset}
\qquad
\deduce[(\Srm^{X}_{T})]{}{}\infer{\Gamma[X/T]}{\Gamma}
$$

\end{enumerate}

\subsection{Embedding function and finite notations for infinitary derivations}

We recall an embedding function $()^{\infty}$ from 
derivations in $\textrm{BI}$ into the derivations in $\textrm{BI}^{\Omega}$, and 
 functions $tp(d)$ and $d[i]$ by Buchholz's method of ``finite notations for infinitary derivations". 
 
The notion of $dg(d)$ where $d \in \textrm{BI}$ is defined as in \cite{Buchholz01} so that $dg(d) \leq dg(d^{\infty})$.

The embedding function $()^{\infty}$ replaces $\R_{C}, \E,\D,\Srm^{X}_{T}$ by $\mathcal{R}_{C}, \mathcal{E}, \mathcal{D},\mathcal{S}^{X}_{T}$ respectively. 
The function $()^{\infty}$ is very similar to the function $()^{\infty}$ in \cite{Buchholz01} 
with the following replacement for induction axiom clause: $(\bigwedge_{\forall x A}(d_{i})_{i \in \omega})^{\infty}: =\bigwedge_{\forall x A}(d_{i}^{\infty}) $.

A derivation $d$ in $\textrm{BI}$ is called \textit{proper} in \cite{Buchholz01} if for every \textit{subderivation} $h$ of $d$,
\begin{enumerate}

\item if $h= \D(h_{0})$, then $dg(h_{0}) = 0$ and $\Gamma(h_{0})$ is an arithmetical sequent;

\item if $h = \Srm^{X}_{T}(h_{0})$, then $h_{0} = \D(h_{00})$.
 
\end{enumerate}
Thus any derivation $d \in \textrm{BI}^{-}$  is proper because it does not contain $\D, \Srm^{X}_{T}$.

Following Buchholz, we define $tp(d)$ and $d[i]$ where $i \in |tp(d)|^{*}$ for each proper derivation $d \in \textrm{BI}$ such that $tp(d)$ is the last inference symbol of $d^{\infty}$, and $d[i]^{\infty}$ is the $i$-th immediate subderivation of $d^{\infty}$.
The definition of $ |tp(d)|^{*}$ is the same as in \cite{Buchholz01}. There is only one new clause in defining $tp(d)$ and $d[i]$:
if $d = (\bigwedge_{\forall x A}(d_{i})_{i \in \omega})$, then $tp(d) := \bigwedge_{\forall x A}, d[i] := d_{i}$ for $i \in \omega$.

\section{Cut-Elimination Theorem for $\textrm{BI}^{-}$}

\subsection{Reduction Relation and Normalization Tree}

 Define 
\[\textrm{BI}_{0} : = \{ \D(d) :    \textrm{where} \ D(d) \ \textrm{is a proper derivation in} \ \textrm{BI} \}.
\]
If $d \in \textrm{BI}_{0}$, then $tp(d) \not \in \{ Cut_{C}, \Omega, \widetilde{\Omega} \}$.

Let $last(d)$ denote the last inference symbol of $d$.
\begin{definition}
$red(d)$
\end{definition}
For a derivation $d \in \textrm{BI}_{0}$ define one step reduction $red(d)$ resulting in a derivation in $\textrm{BI}_{0}$ with $\Gamma(red(d)) = \Gamma(d)$. 

$
red(d)= 
\begin{cases}
\textrm{Ax}_{\Delta} & \textrm{if} \ tp(d) = \textrm{Ax}_{\Delta}; \\
d[0] & \textrm{if} \ tp(d) = \Rep; \\
tp(d)(d[i])_{i \in |tp(d)| } & \textrm{otherwise.} \ 
\end{cases}
$

\begin{definition}
\emph{Let  $ \textrm{BI}^{\Omega} \ni d = I(d_{i})_{i \in |I|}$. Then $|d|$ (the ordinal height of $d$) 
is defined by $|d| := sup(|d_{i}| +1)_{i \in |I|}.$}
\end{definition}

For $d \in \textrm{BI}_{0}$ with $tp(d) = \Rep$, we have  
$|d^{\infty}| > |red(d)^{\infty}|, dg(d) \geq dg(red(d))$, and $\Gamma(d) = \Gamma(red(d))$.

\begin{definition}
Let $d$ be a derivation in $\textrm{BI}_{0}$. We
define the normalization tree $T_{d}^{+}$ as follows. Nodes of the tree are finite sequences $a$ of  natural numbers labeled by derivations $d_{a}$.
\end{definition}

\begin{enumerate}

\item $d^{+}_{\emptyset} : = d$ where $\emptyset$ is the root of $T^{+}_{d}$.

\item If $d^{+}_{a}$ for some $a \in T^{+}_{d}$ is already defined, then the immediate predecessors $d_{a_{i}}^{+}$ of the node $a$ (where $i \in |tp(d_{a}^{+})|$) are defined by cases according to $tp(d_{a}^{+})$.

\begin{enumerate}
\item $tp(d_{a}^{+}) = \Rep$. 

$d_{a_{0}}^{+}: = red(d_{a}^{+})$.

\item $d_{a}^{+}  =  \textrm{Ax}_{\Delta}$. 

In this case, $d_{a}^{+}$ is the leaf of the tree: there is no predecessor of $d_{a}^{+}$.

\item $d_{a}^{+} \neq \textrm{Ax}_{\Delta}$ and $tp(d_{a}^{+}) = \textrm{Ax}_{\Delta}$.

$d_{a_{0}}^{+} : = red(d_{a}^{+})$.

\item Otherwise.

$d_{a_{i}}^{+} : = red(d_{a}^{+})_{i}$ for all $i \in |tp(d_{a}^{+})|$.

\end{enumerate}
\end{enumerate}


\begin{definition}
For every $ d  \in \textrm{BI}_{0}$ let $d^{-}$ be the result of  deleting $\E,\D$ from $d$.
 The normalization tree $T_{d}$ for $d \in \textrm{BI}^{-}$ is the result of
replacing every $d_{a}^{+} \in T_{\D(\E^{m}(d))}^{+}$ by $(d_{a}^{+})^{-}$.
\end{definition}

\begin{definition}
A derivation $d$ of an arithmetical sequent in $\textrm{BI}^{-}$ is normalizable iff $T_{d}$ is well-founded.
\end{definition}

\begin{proposition}
If $T_{d}$ is well-founded, then a cut-free derivation of $\Gamma(d)$ is obtained by deleting some parts of $T_{d}$. 
\end{proposition}
\textbf{Proof}.
By induction on the well-founded normalization tree $T_{d}$.
$\square$

\begin{lemma}
$T^{+}_{d}$ is well-founded for any $d \in \textrm{BI}_{0}$.
\end{lemma}
\textbf{Proof.}
By induction on the ordinal $|(d^{+}_{a})^{\infty}|$. $\square$

\begin{theorem}
Every derivation in $\textrm{BI}^{-}$ of an arithmetical sequent is normalizable.
\end{theorem}
\textbf{Proof.}
It is obvious that $T_{d}$ is well-founded iff $T^{+}_{\D(\E^{m}(d))}$  is well-founded. 
Now apply Lemma 1. $\square$
 
 \vspace{\baselineskip}
 
 From Proposition 1 and Theorem 1, every derivation in $\textrm{BI}^{-}$ of an arithmetical sequent is reduced to a cut-free proof denoted by  $\underline{T_{d}}$.
 For a derivation $d \in \textrm{BI}^{-}$ let $norm(d): = \mathcal{D}(\mathcal{E}^{m}(d))$ for $m = deg(d)$. Then it is easy to see that 
 $\underline{T_{d}} = norm(d)$ up to $\Rep$-inferences.
 
\subsection{Normalization Theorem for Standard Reductions}

The tree $T^{+}_{d}$ is defined as a normalization tree for reductions of derivations in $\textrm{BI}$ containing symbols $\E, \D, \Srm^{X}_{T}$. We analyse what kind of reduction this provides for derivations in $\textrm{BI}^{-}$.

We use Buchholz's notions of {\em nominal form} and $\{ \textrm{R}, \E \}$-form from \cite[pp.~266--267]{Buchholz01}. 
$\mathfrak{a}\{ h \}$ is the result of substitutiong $h$ for $\diamond$ in $\mathfrak{a}:$
\small
\[\deduce{\mathfrak{a}}{\deduce{\vdots}{\dots h \dots}}
\]
\normalsize
Moreover, we adopt his notation $C[k]$ meaning $C_{k}$ for $k \in \{ 0,1 \}$ or $C_{0}(k)$ when $C = C_{0} \wedge C_{1}$ or $\forall x C_{0}(x)$.

Nominal form $\mathfrak{a}\{ h \}$ describes a derivation having a subderivation $h$ and such
that below $h$ only cut-elimination operations $\mathcal{R}_{C}, \mathcal{E}$, the substitution operation
$\mathcal{S}$ and collapsing operation $\mathcal{D}$ are applied. Notation $d=\mathfrak{a}\{ h \}$ indicates that
the end-sequent of $h$ is situated in the end-piece of the derivation $d$. So
for example the condition 2 in the Normalization Theorem below states that the
end-piece of $d_a$ contains a cut suitable for an ``essential $\bigwedge$-reduction" and
this reduction results in the "next" derivation $d_{a_{0}}$. The condition 6
of that theorem states that the end-piece of $d_a$ contains an explicit rule
(cf. the Introduction), and this explicit rule is moved to the bottom of
the derivation. The condition 5 states the cases of axiom-reduction and ``weakening reduction".
\begin{definition}
$FO := \{   \wedge_{A_{0} \bigwedge A_{1}} \bigvee_{A_{0} \vee A_{1}}, \bigwedge_{\forall x  A}, \bigvee_{\exists x A} \}$
\end{definition}

Observe the following simple lemma:

\begin{lemma}
If $d \in \textrm{BI}_{0}$, then $tp(d) \in \{ \textrm{Ax}_{\Delta}, \Rep\} \cup FO$. 
\end{lemma} 

By $d^{+}_{a}$ we denote the derivation attached to a node $a$ in $T_{d}^{+}$. 
Notice that if $tp(d^{+}_{a}) = \textrm{Ax}_{\Delta}$, then there is no successor in $T_{d}^{+}$.

\begin{proposition}
If $tp(d^{+}_{a}) = \Rep$, then one of the following cases holds:
\end{proposition}
\begin{enumerate}

\item $d_{a}^{+} = \mathfrak{a} \{ \E\mathfrak{b} \{ \R_{C}(h_{0},h_{1}) \} \} $, $d_{a_{0}}^{+} = \mathfrak{a} \{ \R_{C[k]}( \E \mathfrak{b} \{  \R_{C}(h_{0}^{-},h_{1}) \} \E \mathfrak{b} \{ \R_{C}(h_{0}, h_{1}^{-}) \} ) \}$.





\item $d_{a}^{+} = \mathfrak{a} \{  \R_{C}(h_{0}, h_{1}) \},$
$
  d^{+}_{a_{0}} = \left\{ \begin{array}{ll}
    \mathfrak{a} \{ h_{i} \}  & \textrm{if} \ tp(h_{1-i}) = \textrm{Ax}_{ \{ C, \neg C \} } \\
    \mathfrak{a} \{ h_{1-i} \}  & otherwise
  \end{array} \right.
  $

\item $d_{a}^{+} = \mathfrak{a} \{ \D h \}$, $d_{a_{0}}^{+} = \mathfrak{a} \{  (\D h)[0]  \}$ with $tp(h) = \widetilde{\Omega}_{\neg \forall X A }$.

\end{enumerate}
\textbf{Proof.} By cases according to  the definition of $tp(d^{+}_{a})$ and $d^{+}_{a}[i]$ as in \cite[p.267]{Buchholz01}. 
$\square$

\begin{proposition}If $I =tp(d^{+}_{a}) \in FO$, then
there is a nominal form $\mathfrak{a}$ such that
\end{proposition}
\begin{enumerate}

\item[] $d^{+}_{a} = \mathfrak{a} \{ {I(h_{i})_{i \in |I|}} \}$, $d^{+}_{a_{i}} = \mathfrak{a} \{ h_{i} \} $ for $i \in |I|$.

\end{enumerate}
\textbf{Proof.} By induction on $d^{+}_{a}$. $\square$

\vspace{\baselineskip}

Since $tp(d^{+}_{a}) \in \{ \textrm{Ax}_{\Delta}, \Rep \} \cup FO $ by Lemma 1, we see that $T_{d}$ describes Gentzen-Takeuti reduction augmented with the reduction of pushing down explicit inferences in the end-piece into the end of derivation.

\begin{theorem}[Normalization by Standard Reductions]
Let $d_{a} \in T_{d}$ for $d \in \textrm{BI}^{-}$. Then $d_{a}$ is an axiom or one of the following cases holds:
\end{theorem}
\begin{enumerate}

\item $d_{a_{0}} = \textrm{Ax}_{\Delta}$.

 

\item $d_{a} = \mathfrak{a} \{  \R_{C}(\mathfrak{b}_{1} \{ \bigwedge_{C_{0} \wedge
  C_{1}}(d_{00}, d_{01}) \}, \mathfrak{b}_{2} \{ \bigvee^{k}_{\neg C_{0} \vee \neg
  C_{1}}(d_{10}) \} ) \}  $, 
 \[d_{a_{0}} = \mathfrak{a} \{ \R_{C[k]}  \mathfrak{b}_{1} \{  \R_{C}(d_{0k},d_{1}) \}  \mathfrak{b}_{2} \{ \R_{C}(d_{0}, d_{10}) \} \}. \]

\item $d_{a} = \mathfrak{a} \{ \R_{C}(\mathfrak{b}_{1} \{ \bigwedge_{\forall x C_{0}(x)}(d_{0n})_{n \in  \omega} \} , \mathfrak{b}_{2} \{ \bigvee^{k}_{\exists x \neg C_{0}(x)}(d_{10}) \} ) \} $, 
 \[d_{a_{0}} = \mathfrak{a} \{ \R_{C[k]}  \mathfrak{b}_{1} \{  \R_{C}(d_{0k},d_{1}) \}  \mathfrak{b}_{2} \{ \R_{C}(d_{0}, d_{10}) \} \}. \]

\item $d_{a} = \mathfrak{a}  \{ \R_{C}(\mathfrak{b}_{1} \{ \bigwedge_{\forall X C_{0}(X)}(d_{00}) \},
\mathfrak{b}_{2} \{ \bigvee^{T}_{\exists X \neg C_{0}(X)}(d_{10}) \}) \}$, 
\[d_{a_{0}} = \mathfrak{a} \{ \R_{C}(\mathfrak{b}_{1} \{ d_{0} \}, \mathfrak{b}_{2} \{ \R_{C_{0}(T)}(\Srm^{X}_{T}(R_{C}(d_{00}, d_{1})), d_{10})) \} ) \}.
\]

\item $d_{a} = \mathfrak{a}\{  \R_{C}(h_{0}, h_{1}) \};$ 
$
  d_{a_{0}} = \left\{ \begin{array}{ll}
     \mathfrak{a} \{ h_{i} \}  & \textrm{if} \ h_{1- i} = \textrm{Ax}_{ \{ C, \neg C \} } \\
     \mathfrak{a} \{ h_{1-i} \}  & otherwise
  \end{array} \right.
  $

\item $d_{a} = \mathfrak{a} \{ {I(h_{i})_{i \in |I|}} \}$, $d_{a_{i}} = \mathfrak{a} \{ h_{i} \} $ for $i \in |I|$.

\end{enumerate}
\textbf{Proof.} By Lemma 2 and Propositions 2,3. $\square$

\vspace{\baselineskip}

As a conclusion we note that derivations in the tree $T^{+}$ satisfy the
additional condition corresponding to Takeuti's requirement in \cite[p.~324, Definition 27.12 (1)]{Takeuti87}) that
all substitution inferences occur in the end-piece. 
\begin{proposition}
If $d_{a} \in T^{+}_{\D(\E^{m}(d))}$ with $d \in \textrm{BI}^{-}$, then all susbstitution inferences in $d_a$
are below all logical inferences, hence the number of substitution rules in $d_{a}$ is finite.  
\end{proposition}
\textbf{Proof.} By bottom-up induction on $T^{+}_{\D(\E^{m}(d))}$. Induction base: the derivation $\D(\E^{m}(d))$ satisfies the condition  since there is no substitution inference in $d$. 
Induction step follows from Propositions 2 and 3.  $\square$

\vspace{\baselineskip}

To illustrate the most important case 4 of cut-elimination in more detail, we present proof figures from \cite[pp.~267--268]{Buchholz01} 
 in a simplified situation using traditional notation. 
Let $d_{a}^{+}$ be of the following form in the traditional notation: 
$$
\infer[\E^{m+1}]{\Gamma}{\infer[\R_{C}]{\Gamma}{\infer[\bigwedge_{\forall X C_{0}(X)}]{\Gamma, \forall X C_{0}(X)}{\deduce{\Gamma, C_{0}(X),\forall X C_{0}(X)}{\vdots}} & \infer[\bigvee^{T}_{\exists X \neg C_{0}(X)}]{\Gamma, \exists X \neg C_{0}(X)}{\deduce{\Gamma, \neg C_{0}(T), \exists X \neg C_{0}(X)}{\vdots}}}}
$$

Therefore $(d_{a}^{+})^{\infty}$ is of the following form:
\small
$$
\infer[\widetilde{\Omega}]{\Gamma}{\infer[\mathcal{R}]{\infer[\mathcal{E}^{m+1}]{\Gamma, C_{0}(X)}{\Gamma, C_{0}(X)}}{\deduce{\Gamma, C_{0}(X),\forall X C_{0}(X)}{\vdots} & \deduce{\Gamma, \exists X \neg C_{0}(X)}{\vdots}} & \infer[\mathcal{E}^{m+1}]{\dots \Gamma, \Delta_{q}, \dots}{\infer[\mathcal{R}]{ \Gamma, \Delta_{q}}{\deduce{\Gamma, \forall X C_{0}(X)}{\vdots} & \infer[\mathcal{R}]{ \Gamma, \Delta_{q}, \exists X \neg C_{0}(X) }{\infer[\mathcal{S}^{X}_{T}]{\Delta_{q}, C_{0}(T)}{\deduce{\Delta_{q}, C_{0}(X)}{\vdots}} & \deduce{\Gamma, \neg C_{0}(T), \exists X \neg C_{0}(X)}{\vdots}}}}}
$$

\normalsize

On the other hand, $d_{a_{0}}^{+} = red(d_{a}^{+})$ is of the following form:

\small

$$
\infer[\E^{m+1}]{\Gamma}{\infer[\R]{\Gamma}{\infer[\bigwedge_{\forall X C_{0}(X)}]{\Gamma, \forall X C_{0}(X)}{\deduce{\Gamma, C_{0}(X),\forall X C_{0}(X)}{\vdots}} & \infer[\R]{\Gamma, \exists X \neg C_{0}(X)}{\infer{\Gamma, C_{0}(T)}{\infer[\D]{\Gamma, C_{0}(X)}{\infer{\infer[\E^{m+1}]{\Gamma, C_{0}(X)}{\Gamma, C_{0}(X)}}{\deduce{\Gamma, C_{0}(X), \forall X C_{0}(X)}{\vdots} & \deduce{\Gamma, \exists X \neg C_{0}(X)}{\vdots}}}} & \deduce{\Gamma,  \neg C_{0}(T), \exists X \neg C_{0}(X)}{\vdots}} }}
$$

\normalsize
Therefore $(d_{a_{0}}^{+})^{\infty}$ is of the following form:
\small
$$
\infer[\mathcal{E}^{m+1}]{\Gamma}{\infer[\mathcal{R}]{\Gamma}{\infer[\bigwedge_{\forall X C_{0}(X)}]{\Gamma, \forall X C_{0}(X)}{\deduce{\Gamma, C_{0}(X),\forall X C_{0}(X)}{\vdots}} & \infer[\mathcal{R}]{\Gamma, \exists X \neg C_{0}(X)}{\infer{\Gamma, C_{0}(T)}{\infer[\mathcal{D}]{\Gamma, C_{0}(X)}{\infer{\infer[\mathcal{E}^{m+1}]{\Gamma, C_{0}(X)}{\Gamma, C_{0}(X)}}{\deduce{\Gamma, C_{0}(X), \forall X C_{0}(X)}{\vdots} & \deduce{\Gamma, \exists X \neg C_{0}(X)}{\vdots}}}} & \deduce{\Gamma,  \neg C_{0}(T), \exists X \neg C_{0}(X)}{\vdots}} }}
$$
\normalsize
In the $\widetilde{\Omega}$-rule at the end of $(d_{a}^{+})^{\infty}$, we have a derivation of $\Gamma, \Delta_{q}$ for each $q \in |\forall X C_{0}(X)|$. Note that $\Gamma$ is an arithmetical sequent, and  $(\D(\E^{m+1}(\R_{C}(d_{00}, d_{1}))), X) \in |\forall X C_{0}(X)|^{*}$. Let's take  $(\mathcal{D}(\mathcal{E}^{m+1}(\mathcal{R}_{C^{i}}(d_{00}^{\infty}, d_{1}^{\infty}))), X)$ as $q$ and $\Gamma$ as $\Delta_{q}$. Then the $q$-th right premise of the $\widetilde{\Omega}$-rule is $(d_{a_{0}}^{+})^{\infty}$, hence
a subderivation of $(d^{+})^{\infty}$. Thus $|(d^{+})^{\infty}| > |(d_{a_{0}}^{+})^{\infty}| $ holds.

\section{Extended $\Omega$-rule}

In this section we extend $\Omega$-rule and prove the cut-elimination theorem for arbitrary (not only arithmetical) sequents.  We define systems $\textrm{BI}^{\Omega^{+}}_{0}$, $\textrm{BI}^{\Omega^{+}}$ in the language $L^{e,i}$, which is obtained from $L$ by adding superscripts $e,i$ to formulas of $L$: if $A \in L$, then   $A^{e} \in L^{e,i}$ and $A^{i} \in L^{e,i}$. 
{\em Sequents} are pairs $\Gamma ; \Delta$ of finite sets of formulas usually written $\Gamma^{e}, \Delta^{i}$. Formulas in $\Gamma^{e}$ are called {\em explicit}, formulas in $\Delta^{i}$ are {\em implicit}.

$\Gamma^{-}$ means the result of deleting all marks $e,i$ occurring in $\Gamma$. Order of (marked) formulas in a sequent is irrelevant. Ordinary arithmetical and second-order rules are modified (in an inessential way) to preserve both $e$ and $i$ marks, so that almost every rule has two versions: for explicit and implicit principal formula. 
In axioms the marks are not important, for example all formulas $A^{e}, \neg A^{e}; A^{e}, \neg A^{i}; A^{i}, \neg A^{e}; A^{i}, \neg A^{i}$ for atomic $A$ are axioms.

\begin{definition}
The systems $\textrm{BI}^{\Omega^{+}}_{0}$ and $\textrm{BI}^{\Omega^{+}}$.
\end{definition}

\begin{enumerate}

\item $\textrm{BI}^{\Omega^{+}}_{0}$ is arithmetic with $\omega$-rule, a rule for second-order universal quantifier, explicit second-order existential quantifier with marks. In what follows, let $\iota \in \{ e,i \}$.

$$
\deduce[(\textrm{Ax}^{\iota}_{\Delta})]{}{}\ \infer{\Delta}{}
$$
\[\ \textrm{where} \ \Delta^{-} = \{ A \} \subseteq \ \textrm{TRUE or } \ \Delta^{-} = \{ C, \neg C \} \ \textrm{for atomic } C\]

$$ 
\deduce[\Huge{(\bigwedge_{A_{0}^{\iota} \wedge A_{1}^{\iota}})}]{}{}\ \infer{A_{0}^{\iota} \wedge A_{1}^{\iota}}{A_{0}^{\iota} & A_{1}^{\iota}}
\qquad
\deduce[\Huge{(\bigvee^{k}_{A_{0}^{\iota} \vee A_{1}^{\iota}})}]{}{}\ \infer[\textrm{where} \ k  \in \{ 0,1 \} ]{A_{0}^{\iota} \vee A_{1}^{\iota}}{A_{k}^{\iota}} 
$$

$$
\deduce[\Huge{(\bigwedge_{\forall x A^{\iota}})}]{}{}\ \infer{\forall x A^{\iota}}{\dots A(x/n)^{\iota} \dots \ \textrm{for all} \ n \in \omega}  
\qquad
\deduce[\Huge{(\bigvee^{k}_{\exists x A^{\iota}})}]{}{}\ \infer[\textrm{where} \ k  \in \omega]{\exists x A^{\iota} }{A(k)^{\iota}} 
$$

$$
\deduce[\Huge{(\bigwedge_{\forall X A^{\iota}}})]{}{}\ \infer[\textrm{where} \ Y \ \textrm{is an eigenvariable}]{\forall X A^{\iota}}{A(X/Y)^{\iota} } 
\qquad 
\deduce[\Huge{(\bigvee^{T}_{\neg \forall X A^{e}})}]{}{}\ \infer{\neg \forall X A^{e} }{\neg A(T)^{e}} 
$$

\item $\textrm{BI}^{\Omega^{+}}$ is obtained by adding the following rules to $\textrm{BI}^{\Omega^{+}}$.
 
$$
\deduce[(\R_{A})]{}{} \ \infer{\emptyset}{A^{i} & \neg A^{i}} 
 \qquad
\deduce[\Huge{(\Omega^{+}_{\neg \forall X A})}]{}{}\ \infer{ \neg  \forall X A^{i}}{\dots \Delta^{\forall X  A^{e}}_{q} \dots \ (q \in |\forall X  A^{e}|)}
$$

$$
\deduce[\Huge{(\widetilde{\Omega}^{+}_{\neg \forall X A})}]{}{}\ \infer[\textrm{where} \ Y \ \textrm{is an eigenvariable}]{\emptyset}{A(Y)^{i} & \dots \Delta^{\forall X  A^{e}}_{q} \dots \  (q \in |\forall X  A^{e}|)} 
$$

with 

\begin{enumerate}

\item $\Delta^{\forall X  A^{e}}_{(d, X)}$ is a sequent of the form $\Gamma^{e},\Pi^{i}$ such that $\Pi^{i}$ is arithmetical  and $\Delta^{\forall X  A^{e}}_{(d, X)} = \Gamma (d)\backslash \ \{  A(X)^{e} \}$,

\item $|\forall X A^{e}| := \{ (d, X ) | \ d \textrm{\ is a cut-free derivation in } \ \textrm{BI}^{\Omega^{+}}_{0},  X \not \in FV(\Delta^{\forall X  A^{e}}_{(d, X)})  \}$. 

\end{enumerate}
\end{enumerate}

\begin{remark} 
\emph{The domain $|\forall X A^{e}|$ of $\Omega^{+}$-rule contains not only cut-free proofs of arithmetical sequents, but cut-free proofs of arbitrary sequents $\Gamma^{e}, \Pi^{i}$ where $\Pi^{i}$ is arithmetical. Moreover a derivation $d$ in this domain may contain $\bigvee_{\neg \forall X A^{e}}$-inferences. }
\end{remark}

\subsection{Cut-Elimination Theorem for $\textrm{BI}^{\Omega^{+}}$}

Cut-degree $dg(d)$ of a derivation $d \in \textrm{BI}^{\Omega^{+}}$ is defined as the least ordinal $\geq rk(C^{i})$ for all implicit formulas $C^{i}$ in $d$. The relation  $d \vdash_{m} \Gamma^{e}, \Pi^{i}$ 
is defined by induction similarly to \cite{Buchholz01}. Only derivations of finite cut-degree are considered below. 

To define one-step reduction $\mathcal{R}_{C}$, collapsing operator $\mathcal{D}$, we need two lemmas.  

\begin{lemma}
There are operators $\mathcal{M}_{i,e}$ and $\mathcal{M}_{e,i}$ changing marks of arithmetical formulas.
\end{lemma}
\begin{enumerate}

\item If $d \vdash_{m} \Gamma, \Lambda^{i}$ where  $\Lambda^{i}$ is an arithmetical, then 
  $\mathcal{M}_{i,e}(d) \vdash_{m} \Gamma, \Lambda^{e}$, and if $d \in \textrm{BI}^{\Omega^{+}}_{0}$ then $\mathcal{M}_{i,e}(d) \in \textrm{BI}^{\Omega^{+}}_{0}$.
  
 \item If $d \vdash_{m} \Gamma, \Lambda^{e}$ where  $\Lambda^{i}$ is  arithmetical, then 
  $\mathcal{M}_{e,i}(d) \vdash_{m} \Gamma, \Lambda^{e}$, and if $d \in \textrm{BI}^{\Omega^{+}}_{0}$ then $\mathcal{M}_{e,i}(d) \in \textrm{BI}^{\Omega^{+}}_{0}$.
  
 \end{enumerate}
\textbf{Proof.} Induction on $d$. 
In both cases $i/e$ and $e/i$ the new derivation is obtained by changing marks of some arithmetical formulas. Axioms stay axioms, for example $A^{i}, \neg A^{i}$ can go to $A^{\tau}, \neg A^{\sigma}$ with any $\tau, \sigma$. All inference rules are preserved. If $d \in \textrm{BI}^{\Omega^{+}}_{0}$ note that arithmetical $A$ cannot be a principal formula of the $\Omega$-rules. $\square$

\begin{theorem}

There is an operator $\mathcal{R}_{C}$ on derivations in $\textrm{BI}^{\Omega^{+}}$ such that 

\end{theorem}

if $d_{0} \vdash_{m} \Gamma^{e}, \Pi^{i}, C^{i}$, \  $d_{1} \vdash_{m} \Gamma^{e}, \Pi^{i}, \neg C^{i}$ and $rk(C^{i}) \leq m$, then $\mathcal{R}_{C}(d_{0}, d_{1}) \vdash_{m} \Gamma^{e}, \Pi^{i}$.

\vspace{\baselineskip}

\textbf{Proof.} 
We consider only cases where the treatment is different from \cite{Buchholz01}.  Let $I_{0}$ and $I_{1}$ be the last inference symbols of $d_{0}$ and $d_{1}$. 

\begin{enumerate}

\item $d_{0} = \textrm{Ax}_{ \{ C^{i}, \neg C^{\tau} \} }$.

Since $C$ is atomic we can assume $\tau = i$ using Lemma 2. Let $\mathcal{R}_{C}(d_{0}, d_{1}):= d_{1}$.

\item $C^{i} = \forall X B(X)^{i}$ is a principal formula in the last inferences of $d_{0}$ and $d_{1}$. Then $d_{0} = \bigwedge_{\forall X B(X)^{i}}(d_{00})$ and $d_{1} = \Omega^{+}(d_{1q})_{q \in |\forall X B(X)^{e}| }$.

Note that the case $d_{1} = \bigvee_{\neg \forall X B(X)^{e}}(d_{0})$ is excluded since $C$ is marked by $i$. 

We have $d_{00} \vdash_{m} \Gamma^{e}, \Pi^{i}, B(Y)^{i}, C^{i}$, and $d_{1q} \vdash_{m} \Gamma^{e}, \Pi^{i}, \Delta^{C^{e}}_{q}, \neg C^{i}$. As in \cite{Buchholz01} $\mathcal{R}_{C}$ is pushed into the premises, and $ \widetilde{\Omega}^{+}$ is introduced: 
\[
\mathcal{R}_{C}(d_{0}, d_{1}) : = \widetilde{\Omega}^{+}(\mathcal{R}_{C}(d_{00}, d_{1}) , \mathcal{R}_{C}(d_{0}, d_{1q}))_{q \in | \forall X B(X)^{e}|}. \ \square
\]

\end{enumerate}

\begin{theorem}

There is an operator $\mathcal{E}$ on derivations in $\textrm{BI}^{ \Omega^{+}}$ such that 

\end{theorem}

\begin{enumerate}

\item if $d \vdash_{m +1} \Gamma^{e}, \Pi^{i}$, then $\mathcal{E}(d) \vdash_{m} \Gamma^{e}, \Pi^{i}$,

\item $\Gamma(d) = \Gamma(\mathcal{E}(d))$.

\end{enumerate}
\textbf{Proof.}
Familiar iteration of $\mathcal{R}_{C}$ (cf. Theorem 1 in \cite{Buchholz01}). $\square$

\vspace{\baselineskip}

A sequent $\Gamma$ is called $\textit{almost explicit}$ if all $i$-marked formulas $A^{i} \in \Gamma$ are arithmetical.
Now we define a collapsing operator $\mathcal{D}$  for arbitrary almost explicit sequent,  which eliminates $\widetilde{\Omega}^{+}$ if $dg(d)=0$. 

\begin{theorem}
There is an operator $\mathcal{D}$ such that 
\end{theorem}

\begin{enumerate}

\item if $\textrm{BI}^{\Omega^{+}} \ni d \vdash_{0} \Gamma$ where $\Gamma$ is almost explicit, then $\textrm{BI}^{\Omega^{+}}_{0} \ni \mathcal{D}(d) \vdash_{0} \Gamma$,

\item $\Gamma(d) = \Gamma(\mathcal{D}(d))$.

\end{enumerate}

\vspace{\baselineskip}

\textbf{Proof.} By induction on $d$.  Since $dg(d) =0$, bottom-up induction on $d$ shows that  all sequents in $d$ are almost explicit and $\Omega^{+}$-rule is not used. $\Gamma \equiv \Gamma^{e},\Pi^{i}$ where $\Pi^{i}$ is an arithmetical sequent since $\Gamma$ is almost explicit.
Let $I$ be the last inference symbol of $d$. We consider only cases where the treatment is different from \cite{Buchholz01}.
\begin{enumerate}

\item $I = \bigvee_{\neg \forall X A^{e}}^{T}$.

In this case $d = \bigvee_{\neg \forall X A^{e}}^{T}(d_{0})$, hence $\Gamma(d_{0}) = \Gamma, \neg A(T)^{e}$ is almost explicit. By IH we have $\textrm{BI}^{\Omega^{+}}_{0} \ni \mathcal{D}(d_{0}) \vdash_{0} \Gamma, \neg A(T)^{e}$. By applying $\bigvee_{\neg \forall X A^{e}}^{T}$ to $\mathcal{D}(d_{0})$, we get the required derivation.

\item $I = \widetilde{\Omega}^{+}$. 

Then $d = \widetilde{\Omega}^{+}(d_{\tau})_{\tau \in \{ 0 \} \cup |\forall X A(X)^{e}|}$.
Now $d_{0} \vdash_{0} \Gamma^{e}, \Pi^{i}, A(Y)^{i}$ and $\textrm{BI}^{\Omega^{+}}_{0} \ni \mathcal{D}(d_{0}) \vdash \Gamma^{e} ,\Pi^{i}, A(Y)^{i}$ by IH with $Y \notin FV(\Gamma(\mathcal{D}(d_{0}))\backslash \{ A(Y)^{i} \})$. Moreover, we get $\textrm{BI}^{\Omega^{+}}_{0} \ni \mathcal{M}(\mathcal{D}(d_{0})) \vdash \Gamma^{e} ,\Pi^{i}, A(Y)^{e}$ by Lemma 2.  Define
$q_{0} := (\mathcal{M}(\mathcal{D}(d_{0})), Y)$, then $q_{0} \in |\forall X A(X)^{e}|$. 
Hence using IH again, we define 
\[
\mathcal{D}(d) : =  \mathcal{D}(d_{q_{0}}) \in \textrm{BI}^{\Omega^{+}}_{0}. \ \square
\]  

\end{enumerate}

\begin{corollary}[Cut-Elimination for $\textrm{BI}^{\Omega^{+}}$]
If $ \textrm{BI}^{\Omega^{+}} \ni d \vdash_{m} \Gamma $ and $\Gamma$ is almost explicit, then $\textrm{BI}^{\Omega^{+}}_{0} \ni \mathcal{D}(\mathcal{E}^{m}(d)) \vdash_{0} \Gamma$.
\end{corollary}

\textbf{Proof.}
By Theorems 3 and 4.
$\square$

\vspace{\baselineskip}
In our marked language, $A^{\tau}[X/T^{\xi}]$ means $(A[X/T])^{\tau}$ where $A$ is a formula of $L$.
Then the substitution operator $\mathcal{S}^{X}_{T}$ is defined by the same induction as in \cite{Buchholz01}.

\begin{lemma}
For every formula $F$ and any $\tau, \sigma \in \{ e, i \}$  
\[\textrm{BI}^{\Omega^{+}} \vdash_{0} F^{\sigma}, \neg F^{\tau}.\] 
\end{lemma}
\textbf{Proof.} By familiar induction on $F$. The only new case is $F \equiv \forall X A$. Then $\neg F \equiv \neg \forall X A$. The rule $\bigwedge_{\forall X A}$ ``splits" $F$, then the rule $\bigvee_{\neg \forall  X A}$ (if $\tau = e$) or $\Omega^{+}_{\neg \forall X A}$ (if $\tau =i$) is applied to ``split" $\neg F$. In the latter case each premise of $\Omega^{+}_{\neg \forall X A}$ is derived a cut:
$$
\infer[Cut]{\Delta, A^{\sigma}}{A^{\sigma}, \neg A^{i} & \Delta, A^{i}}
$$
But this cut is eliminated by applying the operator $\mathcal{E}$. $\square$

\begin{theorem}
There exists an operator $\mathcal{S}^{X}_{T}$ such that

\end{theorem}
if $ \textrm{BI}^{\Omega^{+}}_{0} \ni d \vdash_{0} \Gamma $ and $\Gamma$ is almost explicit, then $\textrm{BI}^{\Omega^{+}} \ni \mathcal{S}^{X}_{T}(d) \vdash_{0} \Gamma[X/T]$.
\vspace{\baselineskip}

\textbf{Proof.}
By familiar induction on $d$. Every rule goes into the same rule, the axioms are treated by Lemma 4.
$\square$

\subsection{Embedding Function and Cut-Elimination Theorem for $\textrm{BI}^{-}$}

In this section we define an embedding function  from 
derivations in $\textrm{BI}^{-}$ into derivations in $\textrm{BI}^{\Omega^{+}}$. 
It is similar to operation $()^{\infty}$ introduced in \cite{Buchholz01} but its inductive definition takes into account marking of formulas in the end-sequent of $d$ as explicit or implicit. As a result, our operation has two arguments: a derivation $d$ and a marking $m$ assigning marks $e,i$ to the formulas in the end-sequent of $\Gamma(d)$. We use notation $d^{ \infty(m)}$ for the resulting derivation of the marked sequent $(\Gamma(d))^{m}$ in $\textrm{BI}^{\Omega^{+}}$. 

In what follows, we assume that $C$ in an axiom $\textrm{Ax}_{ \{  C, \neg C \} }$  of $\textrm{BI}^{-}$ is atomic. A sequent $A, \neg A$ for arbitrary formula $A$ is derived in $\textrm{BI}^{-}$ in a familiar way. 
 
We define $dg(d)$ where $d$ is a derivation in $\textrm{BI}^{-}$ as a good approximation to $dg(d^{\infty(m)})$ so that $dg(d^{\infty(m)}) \leq dg(d)$. 
\begin{definition}
$dg(d)$ 
\end{definition}

Let $d$ be a derivation in $\textrm{BI}^{-}$.

\vspace{\baselineskip}

$
dg(d) :=
\begin{cases}
max(rk(A(T)), dg(d_{0})) & \textrm{if} \ I = \bigvee^{T}_{\neg \forall X  A(X)}; \\
max(rk(C)  , dg(d_{0}), dg(d_{1})) & \textrm{if} \ I = \R_{C} ; \\
sup \{ dg(d_{\tau}) | \tau \in I \} & \textrm{otherwise.}
\end{cases}
$

\vspace{\baselineskip}

When $m$ is a marking of a conclusion $\Gamma$ of some inference rule and $B$ is a side formula of this rule, we denote by $m[B / \tau]$ ($\tau \in \{ e,i \}$) the marking such that  $m[B/\tau](B) = B^{\tau}$ and $m[B/\tau](A) = m(A)$ if $ A \neq B$.
\begin{definition}
Embedding function $()^{\infty(m)}$ from $\textrm{BI}^{-}$ into $\textrm{BI}^{\Omega^{+}}$.
\end{definition}
Let $d$ be a derivation in $\textrm{BI}^{-}$. The function $()^{\infty(m)}$ is defined as follows.
\begin{enumerate}

\item $(\textrm{Ax}_{\Delta})^{\infty(m)} := \textrm{Ax}_{\Delta^{m}}$.


\item $(\bigvee^{T}_{\neg \forall X A}(d_{0}))^{\infty(m)}\\ :=
\begin{cases}
\bigvee_{\neg \forall X A^{e}}(d_{0}^{\infty(m[A(T)/e])}) & \textrm{if} \ m(\neg \forall X A) = \neg \forall X A^{e}; \\
\Omega^{+}_{\neg \forall X A}(\mathcal{R}_{A(T)}(\mathcal{S}^{X}_{T}((\mathcal{M}_{e,i}(d_{q})), d_{0}^{\infty(m[A(T)/i])}))_{q \in |\neg \forall X A^{e}|} & \textrm{otherwise.}
\end{cases}$

\item $(\R_{C}(d_{0},d_{1}))^{\infty(m)} := \mathcal{R}_{C}(d_{0}^{\infty(m[C/i])}, d_{1}^{\infty(m[\neg C/i])})$.

\item Otherwise; 

$I_{\Delta}(d_{i})_{i \in |I|}^{\infty(m)} : = I_{\Delta^{m}}(d_{i}^{\infty(m[A_{i}/ \tau])})$ if $A_{i} \in \Delta_{i}(I)$ and $m(\Delta(I)) = \Delta(I)^{\tau}$.

\end{enumerate}

\begin{theorem}
If $d \in \textrm{BI}^{-}$, then $d^{\infty(m)} \vdash_{dg(d)} \Gamma(d)^{m}$ for any marking function $m$.
\end{theorem}
\textbf{Proof}. By induction on $d$. Let $\Gamma : = \Gamma(d)$. 
We consider only several cases.
\begin{enumerate}

\item $d = \bigwedge_{A_{0} \wedge A_{1}}(d_{0},d_{1})$.

By IH, $d_{0}^{\infty(m[A_{0}/\tau])} \vdash_{dg(d_{0})} \Gamma^{m}, A_{0}^{\tau}$ and $d_{1}^{\infty(m[A_{1} / \tau])} \vdash_{dg(d_{1})} \Gamma^{m},  A^{\tau}_{1}$, hence we obtain the required derivation by applying $ \bigwedge_{A_{0}^{\tau} \wedge A_{1}^{\tau}}$. $\square$

\item $d = (\bigvee_{\neg \forall X A}(d_{0}))$.

\begin{enumerate}

\item $m(\neg \forall X A) = \neg \forall X A^{e}$.

By IH we have $d_{0}^{\infty(m[\neg A(T)/e])} \vdash_{dg(d_{0})} \Gamma^{m}, \neg A(T)^{e}$.
Applying $\bigvee_{\neg \forall X A^{e}}^{T}$, we get the required derivation. Note that $dg(d_{0}) \leq dg(d)$ (cf. Definition 9).

\item Otherwise.

First $d_{q} \vdash_{0} \Pi^{e}, \Lambda^{i}, A(X)^{e}$ where $X$ occurs only in $A(X)^{e}$ for $q \in |\neg  \forall X A^{e}|$. By Lemma 3 and Theorem 5,  we get  $\mathcal{S}^{X}_{T}(\mathcal{M}_{e,i}(d_{q})) \vdash_{0} \Pi^{e}, \Lambda^{i} ,A(T)^{i}$. Since $d_{0}^{\infty(m[A(T)/i])} \vdash_{dg(d_{0})} \Gamma^{m}, \neg A(T)^{i}$ by IH,  we have 
\[\mathcal{R}_{A(T)}(\mathcal{S}^{X}_{T}(\mathcal{M}_{e,i}(d_{q})), d_{0}^{\infty(m[A(T)/i])}) \vdash_{dg(d)} \Gamma^{m},  \Pi^{e}, \Lambda^{i}
 \ \textrm{for all}  \ q \in |\neg  \forall X A^{e}|.\] 
 Thus, by applying $\Omega^{+}_{\neg \forall X A} $, the required derivation is obtained:
 \[\Omega^{+}_{\neg \forall X A}(\mathcal{R}_{A(T)}(\mathcal{S}^{X}_{T}((\mathcal{M}_{e,i}(d_{q})), d_{0}^{\infty(m[A(T)/i])}))_{q \in |\neg \forall X A^{e}|} \vdash_{dg(d)} \Gamma^{m}.
 \]

\end{enumerate}

\item $d = \R_{C}(d_{0}, d_{1})$.

By IH, we have $d_{0}^{\infty(m[C/i])} \vdash_{dg(d_{0})} \Gamma^{m}, C^{i}$ and $d_{1}^{\infty(m[\neg C / i])} \vdash_{dg(d_{1})} \Gamma^{m}, \neg C^{i}$. Thus we get $\mathcal{R}_{C}(d_{0}^{\infty(m[C/i])}, d_{1}^{\infty(m[\neg C/i])}) \vdash_{dg(d)} \Gamma^{m}$ by Theorem 2.

\end{enumerate}

Let $\overrightarrow{e}$ be the marking function assigning $e$ to each formula $A $ in $L$. Moreover, let $d^{*} $ be the result of deleting all marks in sequents and inference rules of $d$. Then we have the following theorem for derivations of  \textit{arbitrary (not only arithmetical) sequents}:

\begin{theorem}[Cut-Elimination Theorem for $\textrm{BI}^{-}$]
If $d \in \textrm{BI}^{-} $, then \\ $\textrm{BI}^{\Omega^{+}}_{0} \ni \mathcal{D}(\mathcal{E}^{n}(d^{\infty^{\overrightarrow{e}}})) \vdash_{0} \Gamma^{\overrightarrow{e}}$ where $n= dg(d)$, hence $\textrm{BI}^{-} \ni (\mathcal{D}(\mathcal{E}^{n}(d^{\infty^{\overrightarrow{e}}})))^{*} \vdash_{0} \Gamma.$
\end{theorem}
\textbf{Proof}.
By Theorem 6 and Corollary 1, we have $\textrm{BI}^{\Omega^{+}}_{0} \ni \mathcal{D}(\mathcal{E}^{n}(d^{\infty^{\overrightarrow{e}}})) \vdash_{0} \Gamma^{\overrightarrow{e}}$.
Since the inference rules in $\textrm{BI}^{\Omega^{+}}_{0}$ become ones of $\textrm{BI}^{-}$ after deleting marks, we get $\textrm{BI}^{-} \ni (\mathcal{D}(\mathcal{E}^{n}(d^{\infty^{\overrightarrow{e}}})))^{*} \vdash_{0} \Gamma$
. $\square$

\bibliographystyle{plain}
\bibliography{mybib}

\begin{thebibliography}{10}

\bibitem{Arai85}
Toshiyasu Arai.
\newblock A subsystem of classical analysis proper to {T}akeuti's reduction
  method for {$\Pi^{1}_{1}$}-analysis.
\newblock {\em Tsukuba Journal of Mathematics}, 9(1):21--29, 1985.

\bibitem{Buchholz81a}
Wilfried Buchholz.
\newblock The {$\Omega_{\mu +1}$}-rule.
\newblock In Wilfried Buchholz, Solomon Feferman, Wolfram Pohlers, and Wilfried
  Sieg, editors, {\em Iterated Inductive Definitions and Subsystems of
  Analysis: Recent Proof-Theoretical Studies}, volume 897 of {\em Lecture Notes
  in Mathematics}, pages 188--233. Springer, 1981.

\bibitem{Buchholz01}
Wilfried Buchholz.
\newblock Explaining the {G}entzen-{T}akeuti reduction steps.
\newblock {\em Archive for Mathematical Logic}, 40:255--272, 2001.

\bibitem{Buchholz_Schutte88}
Wilfried Buchholz and Kurt Sch{\"{u}}tte.
\newblock {\em Proof Theory of Impredicative Subsystems of Analysis}.
\newblock Monographs 2. Bibliopolis, 1988.

\bibitem{Gentzen36}
Gerhard Gentzen.
\newblock Die {W}iederspruchsfreiheit der reinen {Z}ahlentheorie.
\newblock {\em Mathematische Annalen}, 112:494--565, 1936.
\newblock English translation in \cite{Szabo69}.

\bibitem{Gentzen38}
Gerhard Gentzen.
\newblock Neue {F}assung des {W}iderspruchsfreiheitsbeweises f{\"{u}}r die
  reine {Z}ahlentheorie.
\newblock {\em Forschungen zur Logik und zur Grundlegung der exakten
  Wissenschaften}, Neue Folge 4:19--44, 1938.
\newblock English translation in \cite{Szabo69}.

\bibitem{Girard87}
Jean-Yves Girard.
\newblock {\em Proof Theory and Logical Complexity, Vol. {I}}, volume~1 of {\em
  Studies in Proof Theory}.
\newblock Bibliopolis, 1987.

\bibitem{Howard72}
William~A. Howard.
\newblock A system of abstract constructive ordinals.
\newblock {\em Journal of Symbolic Logic}, 37:355--374, 1972.

\bibitem{Mints75a}
Grigori Mints.
\newblock Finite {I}nvestigations of {T}ransfinite {D}erivations.
\newblock 1975.
\newblock In: \cite{Mints92}. Russian original: Zapiski {N}auchnyh {S}eminarov
  LOMI, 49, 67--122.

\bibitem{Mints92}
Grigori Mints.
\newblock {\em {S}elected {P}apers in {P}roof {T}heory}.
\newblock Bibliopolis, 1992.

\bibitem{Pohlers81}
Wolfram Pohlers.
\newblock Proof-theoretical analysis of {$\textrm{ID}_{\nu}$} by the method of
  local predicativity.
\newblock In Wilfried Buchholz, Solomon Feferman, Wolfram Pohlers, and Wilfried
  Sieg, editors, {\em Iterated Inductive Definitions and Subsystems of
  Analysis: Recent Proof-Theoretical Studies}, volume 897 of {\em Lecture Notes
  in Mathematics}, pages 261--357. Springer, 1981.

\bibitem{Schwichtenberg77}
Helmut Schwichtenberg.
\newblock Proof theory: some applications of cut-elimination.
\newblock In Jon Barwise, editor, {\em Handbook of Mathematical Logic},
  volume~90 of {\em Studies in Logic and the Foundations of Mathematics}, pages
  867--895. North-Holland, 1977.

\bibitem{Szabo69}
M.E. Szabo, editor.
\newblock {\em The Collected papers of Gerhard Gentzen}.
\newblock Studies in Logic and the Foundations of Mathematics. North-Holland,
  1969.

\bibitem{Tait65}
William~W. Tait.
\newblock Infinitely long terms of transfinite type.
\newblock In M.~A.~E. Dummett and J.~N. Crossley, editors, {\em Formal Systems
  and Recursive Functions}, Studies in Logic and the Foundations of
  Mathematics, pages 176--185. North-Holland, 1965.

\bibitem{Takeuti67}
Gaisi Takeuti.
\newblock Consistency proofs of subsystems of classical analysis.
\newblock {\em The {A}nnals of {M}athematics}, 86(2):299--348, 1967.

\bibitem{Takeuti87}
Gaisi Takeuti.
\newblock {\em Proof Theory}, volume~81 of {\em Studies in Logic and the
  Foundations of Mathematics}.
\newblock North-Holland, 2nd edition, 1987.

\bibitem{Yasugi70}
Mariko Yasugi.
\newblock Cut elimination theorem for second order arithmetic with the
  {$\Pi^{1}_{1}$}-comprehension axiom and the $\omega$-rule.
\newblock {\em Journal of the Mathematical Society of Japan}, 22(3):308--324,
  1970.

\end{thebibliography}
\end{document}